\documentclass[12pt]{article}
\pdfoutput=1
\usepackage[utf8]{inputenc}
\usepackage{fullpage}
\usepackage[small,bf]{caption}
\usepackage{float}
\usepackage{geometry}
\usepackage{verbatim}
\usepackage{listings}
\usepackage{graphicx}
\usepackage{algorithm}
\usepackage{algpseudocode}
\usepackage{enumerate, enumitem}
\usepackage{booktabs}
\usepackage{hyperref}
\hypersetup{colorlinks=true,
            urlcolor=blue,
            linkcolor=blue,
            citecolor=red}
\usepackage{url}
\usepackage{xcolor}
\usepackage{tikz}
\usetikzlibrary{math, positioning, calc, shapes, arrows.meta, patterns}
\usepackage{xparse}
\usepackage{amsmath,amsfonts}
\usepackage{subcaption}
\usepackage{textcomp}
\usepackage{caption}
\usepackage{threeparttable}
\usepackage{siunitx}
\usepackage{authblk}
\usepackage{multicol}

\usepackage{xargs} 
\usepackage{ifthen} 
\newcommandx{\TODO}[2][1=]{\textcolor{red}{TODO\ifthenelse{\equal{#1}{}}{}{ (#1)}: #2}}
\newcommandx{\NOTE}[2][1=]{\textcolor{blue}{NOTE\ifthenelse{\equal{#1}{}}{}{ (#1)}: #2}}

\definecolor{codegreen}{rgb}{0,0.6,0}
\definecolor{codegray}{rgb}{0.5,0.5,0.5}
\definecolor{codepurple}{rgb}{0.58,0,0.82}
\definecolor{backcolour}{rgb}{0.95,0.95,0.98}

\lstdefinestyle{mystyle}{ 
    backgroundcolor=\color{backcolour},   
    commentstyle=\color{codegreen}, 
    keywordstyle=\color{magenta},
    numberstyle=\tiny\color{codegray}, 
    stringstyle=\color{codepurple},
    basicstyle=\ttfamily\small, 
    breakatwhitespace=false,         
    breaklines=true,                 
    captionpos=t,                    
    keepspaces=true,                 
    numbers=left,                    
    numbersep=5pt,                  
    showspaces=false,                
    showstringspaces=false, 
    showtabs=false,                  
    tabsize=2}

\definecolor{seagreen}{rgb}{0.18, 0.55, 0.34}
\definecolor{mediumviolet-red}{rgb}{0.78, 0.08, 0.52}
\definecolor{khaki}{rgb}{0.94, 0.9, 0.55}

\lstdefinelanguage{mypython}
{
	keywords=[1]{from, import, as, assert, not, print, nonneg, boolean},
	keywordstyle=[1]{\color{mediumviolet-red}},
	keywords=[2]{cp, np, pd, cvxpy, numpy, pandas, sum, multiply, Variable, array,
    max, unique, sum_largest, hstack, matmul, Problem, Minimize, solve, range, set,
    max_daily_powers},
	keywordstyle=[2]{\color{seagreen}},
	upquote=true,
	showstringspaces=false,
	basicstyle=\ttfamily,
	columns=fullflexible,
	keepspaces=true,
	emph={True,False,def,return,float,class,match,switch,len},
	emphstyle={\color{seagreen}},
	belowskip=1em,
	aboveskip=1em,
    morecomment=[l]{\#}
}

\usepackage[
backend=bibtex,
style=numeric-comp,
bibstyle=ieee,
sorting=ynt,
firstinits=true,
url=false, doi=false, eprint=false,
]{biblatex}
\newcommand{\ie}{{\it i.e.}}

\newcommand{\BEQ}{\begin{equation}}
\newcommand{\EEQ}{\end{equation}}
\newcommand{\BEAS}{\begin{eqnarray*}}
\newcommand{\EEAS}{\end{eqnarray*}}
\newcommand{\ones}{\mathbf 1}

\newcommand{\reals}{{\mbox{\bf R}}}

\newcommand{\naturals}{{\mbox{\bf N}}}

\newcommand{\symm}{{\mbox{\bf S}}}  

\newcommand{\expect}[1]{\mathop{\bf E{}}[#1]}

\newcommand{\argmin}{\mathop{\rm argmin}}

\usepackage{relsize} 

\newcommand{\BIT}{\begin{itemize}}
\newcommand{\EIT}{\end{itemize}}
\newcommand{\BNUM}{\begin{enumerate}}
\newcommand{\ENUM}{\end{enumerate}}

\begin{document}

\title{An Operator Splitting Method for Large-Scale CVaR-Constrained Quadratic Programs}

\author[1]{Eric Luxenberg\footnote{Equal contribution.}}
\newcommand\CoAuthorMark{\footnotemark[\arabic{footnote}]} 
\author[1,2]{David Pérez-Piñeiro\protect\CoAuthorMark}
\author[1]{Steven Diamond}
\author[1]{Stephen Boyd}
\affil[1]{Department of Electrical Engineering, Stanford University}
\affil[2]{Department of Chemical Engineering, Norwegian University of Science and Technology}

\maketitle

\begin{abstract} 
We introduce a fast and scalable method for solving quadratic programs with
conditional value-at-risk (CVaR) constraints. While these problems can be
formulated as standard quadratic programs, the number of variables and
constraints grows linearly with the number of scenarios, making general-purpose
solvers impractical for large-scale problems. Our method combines operator
splitting with a specialized $O(m\log m)$ algorithm for projecting onto CVaR
constraints, where $m$ is the number of scenarios. The method alternates between
solving a linear system and performing parallel projections, onto CVaR
constraints using our specialized algorithm and onto box constraints by
simple clipping. Numerical examples from several application domains
demonstrate that our method outperforms general-purpose solvers by several
orders of magnitude on problems with up to millions of scenarios. Our method is
implemented in an open-source package called CVQP.
\end{abstract}

\clearpage
\tableofcontents
\clearpage

\section{Introduction}
Many applications in finance and engineering require controlling the risk
of extreme outcomes. A widely used measure for tail risk is the
\emph{conditional value-at-risk} (CVaR), defined as the expected value of
losses exceeding a given quantile. CVaR is a coherent and convex risk
measure~\cite{rockafellar2000optimization}, so optimization
problems involving CVaR can be reliably and efficiently solved.

Many practical applications, from portfolio optimization to quantile regression,
can be formulated as quadratic programs with CVaR constraints. While these
problems are convex and can be reformulated as standard quadratic programs, the
number of variables and constraints grows linearly with the number of scenarios.
For problems with many scenarios, general-purpose solvers become prohibitively
slow or fail entirely. To address this challenge, we develop a fast and scalable
method for solving quadratic programs with CVaR constraints.

We present two main contributions. First, we develop an $O(m\log m)$ algorithm
for projecting onto CVaR constraints, where $m$ is the number of scenarios.
Building on this algorithm, our second contribution is an operator splitting
method for solving large-scale CVaR-constrained quadratic programs. The method
alternates between solving a linear system and performing parallel projections,
onto CVaR constraints using our specialized algorithm and onto box constraints
by simple clipping. Numerical examples from
several application domains demonstrate that our method outperforms
general-purpose solvers by several orders of magnitude on problems with up to
millions of scenarios.

\subsection{Conditional value-at-risk (CVaR)}
The conditional value-at-risk (CVaR) at level $\beta$ is a risk measure that
captures the expected value over the worst $(1-\beta)$ fraction of outcomes of a
real-valued random variable. For a random variable $X$ representing losses
(where larger values are worse), we first define the value-at-risk (VaR) at
level $\beta$ as
\[
    \psi_\beta(X) = \inf\{x\mid \mathbb{P}(X\leq x) \geq \beta\}.
\]
CVaR at level $\beta$ is defined as the expected value of all losses exceeding
VaR at level $\beta$, 
\[
    \phi_\beta(X) = \expect{X\mid X\geq \psi_\beta(X)}.
\]
At the same level $\beta$, CVaR provides an upper bound on VaR.

\paragraph{Sample CVaR.}
For a finite set of samples $z_1,\ldots, z_m\in\reals$ representing the distribution
of losses, Rockafellar and Uryasev \cite{rockafellar2000optimization} showed that
the CVaR at level $\beta$, denoted $\phi_\beta:\reals^m \to \reals$, can be
computed as
\[
    \phi_\beta(z) = \inf_{\alpha \in \reals} \left\{ \alpha + \frac{1}{(1-\beta)m}\sum_{i=1}^m(z_i-\alpha)_+ \right\},
\]
where $(z-\alpha)_+ = \max\{z-\alpha, 0\}$ is the positive part of $z-\alpha$. 
A CVaR constraint at level $\beta$ can be represented by a set of linear
inequality constraints,
\[
    \left\{z\mid \phi_\beta(z) \leq \kappa\right\} = \left\{z ~\middle|~ 
    \begin{array}{l}\alpha + \frac{1}{(1-\beta)m}\sum_{i=1}^m y_i \leq \kappa,\\
    z_i - \alpha \leq y_i,~0 \leq y_i, \quad i = 1,\ldots, m
    \end{array}
    \right\}.
\]

\subsection{CVaR-constrained quadratic programs}
This work presents a customized solver for the \textit{CVaR-constrained
quadratic program} (CVQP),
\BEQ \label{e-cvar_constrained_qp}
\begin{array}{ll}
\mbox{minimize} & (1/2)x^T Px + q^T x \\
\mbox{subject to} & \phi_\beta(Ax) \leq \kappa, \quad l\leq Bx\leq u,
\end{array}
\EEQ 
where $x \in \reals^n$ is the decision variable. The objective function is
defined by a positive semidefinite matrix $P \in \symm^n_{+}$ and a vector $q
\in \reals^n$. The CVaR constraint is defined by
quantile level $\beta \in (0, 1)$, threshold $\kappa \in \reals$, and matrix
$A \in \reals^{m \times n}$. Additional linear inequality constraints are defined by a matrix
$B \in \reals^{p \times n}$, and vectors
$l\in\left({\reals\cup\{-\infty\}}\right)^p$,
$u\in\left({\reals\cup\{\infty\}}\right)^p$. 
Since $P$ is positive semidefinite and the constraints are convex, the CVQP is
a convex optimization problem.
We focus on the case where $p\ll m$,
\ie, the number of additional constraints is much smaller than the number of
scenarios.

We assume the CVQP problem \eqref{e-cvar_constrained_qp} is feasible and
that $P$, $A$, and $B$ have zero common nullspace, which ensures that the
problem is bounded. This condition also implies that the matrix
$M = P + \rho(A^TA + B^TB)$ is positive definite for any $\rho > 0$, which
is needed for the ADMM linear system solve described in \S\ref{s-admm}.
In practice the condition holds when $P$ is positive definite or when $A$
has full column rank.

\paragraph{CVaR terms in the objective.} 
The CVQP formulation extends to problems with CVaR terms in the objective.
Consider the problem
\[
\begin{array}{ll}
\mbox{minimize} & \frac{1}{2}x^T P x + q^T x + \phi_\beta(Ax) \\
\mbox{subject to} & l \leq Bx \leq u.
\end{array}
\]
Using the translation-equivariance property of CVaR (i.e., $\phi_\beta(X + c) = \phi_\beta(X) +
c$), we introduce an auxiliary variable $t$ and reformulate as
\[
\begin{array}{ll}
\mbox{minimize} & \frac{1}{2}x^T P x + q^T x + t \\
\mbox{subject to} & \phi_\beta(Ax-t) \leq 0, \quad l \leq Bx \leq u.
\end{array}
\]
Defining the decision variable as the pair $(x, t)$, this is an instance
of the standard CVQP form \eqref{e-cvar_constrained_qp}.

\paragraph{Quadratic program solvers.} By using the linear inequality
representation of the CVaR constraint, the CVQP can be reformulated as a
standard quadratic program and solved with any QP solver. However, when the
number of scenarios $m$ is large, this approach becomes prohibitively slow.

\paragraph{Applications.}
CVaR-constrained quadratic programs arise in many domains.
In portfolio optimization, $x$ represents portfolio weights and $Ax$
gives portfolio losses across $m$ return scenarios; the CVaR constraint
limits the expected loss in the worst $(1-\beta)$ fraction of scenarios,
while the quadratic objective captures risk-adjusted returns
\cite{krokhmal2002portfolio}.
In supply chain planning, $x$ includes order quantities and routing
decisions, with scenarios modeling demand uncertainty and supplier
disruptions; CVaR constraints bound extreme total cost realizations
\cite{sawik2011selection, gotoh2007newsvendor}.
Similar formulations appear in network flow problems with stochastic
capacities \cite{sorokin2013computational}, facility location and
disaster response logistics \cite{noyan2012risk}, energy systems
with uncertain generation and prices \cite{doege2006risk}, and radiation
treatment planning, where CVaR constraints replace intractable
dose-volume constraints \cite{romeijn2003novel}.
In control and statistical learning, CVaR provides tractable convex
approximations to chance constraints
\cite{van2015distributionally, laguel2021superquantiles}.
For a comprehensive review, see \cite{filippi2020conditional}.

\subsection{Related work}
CVaR constraints were introduced by Rockafellar and Uryasev
\cite{rockafellar2000optimization} and are now standard in risk-constrained
optimization \cite{krokhmal2002portfolio}. Our method builds on operator
splitting \cite{boyd2011distributed}, which decomposes the problem into simpler
subproblems: a linear system solve and projections.

The CVaR projection algorithm adds to a family of finite-termination projection
methods, including the classical simplex projection \cite{gafni1984two} and
isotonic regression \cite{barlow1972isotonic}. Concurrently with the
development of our CVaR projection algorithm, Roth and Cui released a preprint \cite{roth2023n}, 
later published as \cite{roth2025n}, containing two finite termination algorithms for CVaR projection attaining the
same complexity. After corresponding with the authors, we were
unable to find any equivalence between their algorithms and ours.

Roth and Cui \cite{roth2024fast} also proposed a second-order computational
framework for CVaR-constrained optimization, using a semismooth Newton-based
augmented Lagrangian method. As a second-order method, it can achieve superlinear convergence near a
solution, while ADMM is a first-order method with an $O(1/k)$ convergence
rate. On the other hand, each ADMM
iteration is cheap (a single backsolve and parallel projections after an
initial factorization), and the method is simple to implement. In our
experience, ADMM reaches moderate accuracy in few iterations, which is
sufficient for most applications.

\subsection{Outline}
This paper is organized as follows. In \S\ref{s-admm}, we present an ADMM-based
solution method for the CVQP. \S\ref{s-cvar_projection} establishes theoretical
foundations for projecting onto CVaR constraints, which we use to develop our
$O(m\log m)$ projection algorithm in \S\ref{s-procedure}. In
\S\ref{s-experiments}, we present benchmark results comparing our method against
state-of-the-art solvers on portfolio optimization and quantile regression
problems with up to millions of scenarios. \S\ref{s-extensions} discusses
possible extensions to our method.

\section{Solution via ADMM}
\label{s-admm}

\subsection{ADMM}
We solve the CVQP by reformulating the problem
and applying the alternating direction method of multipliers (ADMM).
By introducing auxiliary variables $z \in \reals^m$ and
$\tilde{z} \in \reals^p$, we can rewrite problem
\eqref{e-cvar_constrained_qp} as
\[
\begin{array}{ll}
\mbox{minimize} & (1/2)x^T Px + q^T x + I_{\mathcal{C}}(z) + I_{[l,u]}(\tilde{z}) \\
\mbox{subject to} & Ax = z,\quad Bx = \tilde{z},
\end{array}
\]
where $\mathcal{C} = \{z \mid \phi_\beta(z) \leq \kappa\}$.
The auxiliary variables $z$ and $\tilde z$ decouple the CVaR constraint
from the box constraints, so that each can be handled separately
in the ADMM updates below.
Here, $I_{\mathcal{C}}$ and $I_{[l,u]}$ are the indicator functions for the sets
$\mathcal{C}$ and $\{\tilde z \mid l \leq \tilde z \leq u\}$, respectively, given by 
\[
    I_{\mathcal{C}}(z) = \begin{cases} 0 & z\in\mathcal{C} \\ \infty & \text{otherwise} \end{cases}, 
    \quad I_{[l,u]}(\tilde{z}) = \begin{cases} 0 & l\leq \tilde z \leq u\\ \infty & \text{otherwise} \end{cases}.
\]

We use the scaled dual variables $u=(1/\rho)y$
and $\tilde{u}=(1/\rho)\tilde{y}$, where $\rho > 0$ is a hyperparameter and $y
\in \reals^m$ and $\tilde{y} \in \reals^p$ are the vectors of dual variables
associated with the equality constraints $Ax = z$ and $Bx = \tilde{z}$,
respectively.
We use ADMM with over-relaxation 
(see \cite{eckstein1992douglas, eckstein1994parallel} for
analysis), with over-relaxation parameter $\alpha\in(0,2)$.

For notational convenience we define the linear system parameters
\[
    M = P + \rho (A^TA + B^TB), \quad p^k = q - \rho A^T(z^k - u^k) - \rho B^T(\tilde{z}^k - \tilde{u}^k).
\]
Our assumption on the common nullspace of $P$, $A$, and $B$ implies that
$M$ is positive definite.
The ADMM updates are
\begin{eqnarray}
    x^{k+1} &:=& \argmin_x \left( \frac{1}{2}x^TMx + (p^k)^Tx\right) 
\label{e-update_x-final} \\
    z^{k+1/2} &:=& \alpha A x^{k+1} + (1-\alpha)z^k\label{e-update_z_half}\\
    \tilde{z}^{k+1/2} &:=& \alpha B x^{k+1} + (1-\alpha)\tilde{z}^k\label{e-update_ztilde_half}\\
    z^{k+1} &:=& \Pi_{\mathcal{C}}\left(z^{k+1/2}+ u^k\right) \label{e-update_z_proj-final} \\
    \tilde{z}^{k+1} &:=& \Pi_{[l,u]}\left(\tilde{z}^{k+1/2} + \tilde{u}^k\right)
    \label{e-update_ztilde_proj-final}\\
    u^{k+1} &:=& u^k + z^{k+1/2} - z^{k+1}\label{e-update_u}\\
    \tilde{u}^{k+1} &:=& \tilde{u}^k + \tilde z^{k+1/2} - \tilde{z}^{k+1}.\label{e-update_utilde}
\end{eqnarray}
The operator $\Pi_\mathcal{C}$ is the (Euclidean) projection onto
$\mathcal{C}$ and
$\Pi_{[l,u]}$ is the projection onto the set $[l,u]\subset \reals^p$.
The over-relaxation, dual updates, and box projection are all trivial
operations; the linear system solve requires one cached factorization
of $M$. The main computational bottleneck is the CVaR projection
\eqref{e-update_z_proj-final}, which we address in \S\ref{s-cvar_projection}--\ref{s-procedure}.

It is well known that this ADMM algorithm converges
\cite{eckstein1992douglas,GABAY1983299}, with a worst-case $O(1/k)$
convergence rate, though linear convergence is commonly observed in
practice for quadratic programs \cite{boyd2011distributed}.

\paragraph{Dynamic penalty parameter.} 
While ADMM converges for any fixed penalty parameter $\rho > 0$, the convergence
rate can be improved by adaptively updating $\rho$. A simple and effective
update scheme adjusts $\rho$ based on the relative magnitudes of the primal and
dual residuals. When the primal residual is $\mu$ times larger than the dual
residual, we multiply $\rho$ by a factor $\tau > 1$; when the dual residual is
$\mu$ times larger than the primal residual, we divide $\rho$ by $\tau$. Since
$\rho$ appears in the matrix $M$, each update requires re-factorizing the linear
system. To balance computational cost with convergence benefits, we perform
these updates at fixed intervals of $T$ iterations.

\subsection{Evaluating the updates}
The updates
\eqref{e-update_z_half}, \eqref{e-update_ztilde_half}, \eqref{e-update_u}, and
\eqref{e-update_utilde} are trivial to implement.
The update \eqref{e-update_ztilde_proj-final} is also a very simple clipping
operation:
\[
\tilde z^{k+1}_i = \begin{cases} 
l_i & v_i < l_i \\ 
v_i & l_i \leq v_i \leq u_i \\ 
u_i & v_i > u_i 
\end{cases}
\]
where $v = \tilde z^{k+1/2}+\tilde u^k$.

The update \eqref{e-update_x-final} can be expressed as $x^{k+1} = -M^{-1}p^k$,
which requires the solution of a linear system of equations with positive definite 
coefficient matrix $M$.
We factorize $M$ once
(for example, a sparse Cholesky or $LDL^T$ factorization) so that subsequent solves
require only the backsolve.
This can substantially improve the efficiency of this step; when $M$ is dense,
for example, the first factorization costs $O(n^3)$ flops, while subsequent solves
require only $O(n^2)$ flops.  (When $M$ is sparse the speedup from caching 
the factorization is smaller than $n$, but still very significant.)

That leaves only the CVaR projection \eqref{e-update_z_proj-final}, for
which we develop an $O(m\log m)$ algorithm in the next two sections.

\subsection{Summary}
For convenience we summarize the complete CVQP method in
Algorithm~\ref{a-cvqp}. The CVaR projection is the main computational
bottleneck and uses the $O(m\log m)$ algorithm described in
\S\ref{s-procedure}.

\begin{algorithm}[h]
\caption{CVQP solver}\label{a-cvqp}
\begin{algorithmic}[1]
\State Factorize $M = P + \rho(A^TA + B^TB)$.
\Repeat
    \State Solve linear system \eqref{e-update_x-final}.
    \State Over-relax \eqref{e-update_z_half}--\eqref{e-update_ztilde_half}.
    \State Project onto CVaR constraint \eqref{e-update_z_proj-final}; see \S\ref{s-procedure}.
    \State Clip to box constraints \eqref{e-update_ztilde_proj-final}.
    \State Update dual variables \eqref{e-update_u}--\eqref{e-update_utilde}.
    \State Update $\rho$ and re-factorize $M$ if needed.
\Until{primal and dual residuals below tolerance.}
\end{algorithmic}
\end{algorithm}

\section{CVaR projection preliminaries} \label{s-cvar_projection}
We now present preliminaries that motivate our
efficient algorithm for evaluating the projection operator $\Pi_\mathcal{C}$, \ie, 
for solving the problem
\begin{equation} \label{e-cvar_projection}
\begin{array}{ll}
\mbox{minimize} & \|v - z\|_2^2 \\
\mbox{subject to} & \phi_\beta(z) \leq \kappa,
\end{array}
\end{equation}
with variable $z \in \reals^m$.

The CVaR constraint can be expressed equivalently using the \textit{sum of $k$
largest components}
function, 
\BEQ\label{e-sumkproj}
  f_k(z)=\sum_{i=1}^k z_{[i]}, 
\EEQ
where $z_{[1]}\geq z_{[2]}\geq \cdots \geq z_{[m]}$ are the components of $z$ in 
nonincreasing order. This is a convex function that sums the $k$ largest elements of a vector
(the pointwise maximum of $\binom{m}{k}$ linear functions).
Setting $k = (1-\beta)m$ and evaluating the infimum over $\alpha$ in the
Rockafellar-Uryasev formula (the minimizer is $\alpha = z_{[k]}$), we obtain
$\phi_\beta(z) = (1/k) f_k(z)$. The constraint $\phi_\beta(z) \leq \kappa$ is
therefore equivalent to $f_k(z) \leq \kappa k$.
We rewrite problem \eqref{e-cvar_projection} as
\begin{equation} \label{e-sumklargest_projection}
\begin{array}{ll}
\mbox{minimize} & \|v - z\|_2^2 \\
\mbox{subject to} & f_k(z) \leq d,
\end{array}
\end{equation}
where $d = \kappa k$. We assume $(1-\beta) m\in\naturals$
(and round up to the nearest integer if not).

We will use the formulation \eqref{e-sumklargest_projection} to derive
our CVaR projection algorithm. While this is a computationally tractable convex
optimization problem, we seek closed form or computationally efficient ways to
evaluate this projection operator. Since merely evaluating $f_k$ with $k=(1-\beta) m$ has a
complexity of $O(m \log m)$, we desire a similar complexity for the projection
operator. 

\subsection{Sorted projection}\label{s-perm}
Let $P\in \reals^{m\times m}$ be a permutation matrix that sorts $v$ in
descending order, \ie, $(Pv)_1 \geq (Pv)_2 \geq \cdots \geq (Pv)_m$.
We will denote the sorted vector $Pv$ as $v'$.
Since for any $z\in \reals^m$,
\[
    \|v-z\|_2^2 = \|Pv-Pz\|_2^2, \qquad f_k(z) = f_k(Pz), 
\]
the solution to the projection problem is given by
\BEAS
    z^\star &=& \argmin_{f_k(z) \leq d} \|v-z\|_2^2\\
            &=& \argmin_{f_k(Pz) \leq d} \|Pv-Pz\|_2^2\\
            &=& P^T\left(\argmin_{f_k(u) \leq d} \|v'-u\|_2^2\right).\\
\EEAS
The last equality is the change of variable $u = Pz$.
In words, we can sort $v$ in descending order, project the sorted vector $v'$, and unsort
the result to obtain the projection of $v$ onto the set $\{z\mid f_k(z)\leq
d\}$.

\subsection{Optimality conditions}\label{s-kkt}
We characterize the optimality conditions of the projection problem; these
conditions will be used in \S\ref{s-correctness} to prove the correctness of the
algorithm. Let $\mathcal{A}= \{a\in \{0,1\}^m \mid \ones^Ta = k\}$
be the set of all $\binom{m}{k}$ $k$-hot indicator vectors.
Since $f_k(z) = \max_{a \in \mathcal{A}} a^Tz$, the constraint
$f_k(z) \leq d$ is equivalent to $a^Tz \leq d$ for all $a \in \mathcal{A}$.
We can therefore rewrite problem \eqref{e-sumklargest_projection}
in the equivalent form
\[
\begin{array}{ll}
\mbox{minimize} & \frac{1}{2}\|v - z\|_2^2 \\
\mbox{subject to} & a_i^Tz \leq d, \quad i=1,\ldots,\binom{m}{k},
\end{array}
\]
where the subscript $i$ enumerates the elements of $\mathcal{A}$.

The KKT conditions tell us a $z, \lambda$ pair is optimal if and only if 
\[
    v-z = \sum_{i=1}^{\binom{m}{k}} \lambda_i a_i,
\] 
and
\[
    \qquad a_i^Tz \leq d, \quad \lambda_i \geq 0, 
    \quad \left(a_i^Tz < d \implies \lambda_i = 0\right), \qquad i=1,\ldots,\binom{m}{k}.
\]

The optimality conditions motivate an algorithm that solves the KKT
system. We will define $z = v-\Delta$, where $\Delta$ is the vector we remove
from $v$ to get the projection. The KKT conditions above imply that
$\Delta$ is a nonnegative combination of the vectors $a_i$ for which
$a_i^Tz=d$, \ie, the indicators of the largest $k$ entries of $z$. (Note
that there are potentially several such vectors in the presence of ties: consider the
largest 2 entries in the vector $(2,1,1,0)$.)
Thus, if we can construct a vector $\Delta$ that is a nonnegative combination of
largest-$k$ indicator vectors, such that $f_k(v-\Delta) = d$, 
then we have found the projection $z= v-\Delta$.

\subsection{Algorithm motivation}
We describe at a high level the motivation behind our projection algorithm.
We pre-process by sorting and work with the sorted vector $v'$, since from
\S\ref{s-perm} we know that sort-project-unsort is the same as projecting. We then proceed
by incrementally constructing $\Delta$ until the sum of the $k$ largest elements
of $v'-\Delta$ is equal to $d$.

We want to begin decreasing elements of $v'$ until the sum of the $k$ largest
elements is equal to $d$. Intuitively, we should decrease the largest $k$
elements without modifying the remaining elements, since they do not affect the sum.
We begin by reducing the largest elements uniformly until either the constraint is
satisfied, or until the $k$th largest element ties the next element. 
We reduce each of the $k$ largest elements by the same amount because we are
penalized by the sum of squared changes.

If reducing the largest $k$ elements causes elements $k$ and $k+1$ to tie, it is
intuitive that the tied elements should be further reduced by the same amount;
otherwise a reduction would be applied to an element not in the largest $k$. 
Since we also want to uniformly decrease the largest elements preceding the tie,
we need to decide the ratio of these two reductions. Given this ratio, we reduce
the untied and tied entries until either the constraint is satisfied or the next
tie occurs.

We will need to determine this ratio at each step of the algorithm, since
at a given step of the algorithm, the altered vector $v'-\Delta$ will be partitioned
into three regions: first the largest $n_u\in\naturals$ \emph{untied} elements, then the
$n_t\in\naturals$ \emph{tied} elements, and lastly the remaining unaltered elements:
\[
    v'-\Delta = \left(\underbrace{v'_1 \geq \cdots \geq v'_{n_u}}_{\text{untied}} 
    ~\geq~\underbrace{v'_{n_u+1}=\cdots=v'_{n_u+n_t}}_{\text{tied}}
    ~\geq~\underbrace{v'_{n_u+n_t+1}\geq \cdots\geq v'_m}_{\text{unaltered}}\right).
\]
The core insight of the algorithm is deriving the ratio of reducing the untied
and tied elements; this ratio is derived in \S\ref{s-procedure}.

\section{CVaR projection algorithm}\label{s-procedure}
With the preliminaries in place, we now present our CVaR projection algorithm.
\paragraph{Pre- and post-processing.}
The CVaR projection algorithm is applied to the vector $v$ sorted in descending
order as described in \S\ref{s-perm}. The complexity of sorting $v$ and then unsorting
the projection is $O(m\log m)$. The sorting and unsorting are in practice carried
out with indices; the permutation matrix $P$ is never explicitly formed.

\paragraph{Algorithm state.}
The algorithm's state is
\begin{multicols}{2}
\begin{itemize}
    \item $j$: the iterate number,
    \item $n_u$: the number of untied entries,
    \item $n_t$: the number of tied entries,
    \item $\eta$: the decrease to the untied entries,
    \item $S$: the sum of the largest $k$ entries,
    \item $a_t$: the value of the tied entries,
    \item $a_u$: the value of the last untied entry, 
    \item $a_e$: the value of the first unaltered entry.
\end{itemize}
\end{multicols}
The algorithm is initialized with $n_u=k$, $n_t=0$, $\eta=0$,
$S=\sum_{i=1}^k v'_i$, $a_u=v'_k$, and $a_e=v'_{k+1}$.
The algorithm repeats the \textit{decrease step} described in \S\ref{s-decrease}
until $S = d$, after which the projection of the sorted $v'$,
\[
    z = \left(
        \underbrace{v'_1-\eta, \ldots, v'_{n_u}-\eta}_{n_u}, ~\underbrace{a_t, \ldots, a_t}_{n_t}, ~v'_{n_u+n_t+1}, \ldots, v'_m
    \right)
\] 
is formed, unsorted, and returned.

\subsection{Decrease step}\label{s-decrease}
At each step, the algorithm decreases the entries of $z$ along a direction
$\delta \in \reals^m$ that depends on the current group structure.
The update is $z \leftarrow z - s_0\delta$, where $s_0 > 0$ is chosen
as the largest scaling such that $z$ remains sorted and the sum of its
$k$ largest entries stays at least $d$.
Three events can limit $s_0$:
\BIT
\item the last untied entry reaches the tied value
      (\emph{merge}: the untied block shrinks by one),
\item the tied value reaches the first unaltered entry
      (\emph{absorb}: the tied block grows by one),
\item the sum constraint $S = d$ is met (\emph{termination}).
\EIT
On every nonterminal step, the tied block grows by one entry
(from a merge or an absorb).
The algorithm terminates in at most $m$ steps.
We handle the first step separately because there is no tied block yet;
all subsequent steps share the same general formulas.

\paragraph{First step.}
On the first step, $n_u = k$ and $n_t = 0$.
The direction is
\[
\delta = \left(\underbrace{1, \ldots, 1}_{k},~
\underbrace{0,\ldots,0}_{m-k}\right),
\]
\ie, all $k$ entries decrease uniformly.
Since there is no tied block, only absorb and termination can occur, giving
$s_0 = \min(a_u - a_e,~ (S-d)/k)$.
If $S - s_0 k = d$, the algorithm terminates with $\eta = s_0$.
Otherwise, entries $k$ and $k+1$ form a tied block, and the state
is updated to
\[
\eta \leftarrow s_0, \quad
S \leftarrow S - s_0 k, \quad
n_u \leftarrow k-1, \quad
n_t \leftarrow 2, \quad
a_t \leftarrow v'_{k+1}.
\]
The values of $a_u$ and $a_e$ are updated as in the general step below,
using the new $n_u$, $n_t$, and $\eta$.

\paragraph{General step.}
For all subsequent steps, $n_t \geq 2$ and $k - n_u \geq 1$.
The direction is
\[
\delta(n_u,n_t,k) = \left(\underbrace{\frac{n_t}{k-n_u}, \ldots,
\frac{n_t}{k-n_u}}_{n_u}, ~\underbrace{1, \ldots, 1}_{n_t},~
\underbrace{0,\ldots,0}_{m-n_u-n_t}\right).
\]
The untied entries decrease at rate $n_t/(k-n_u) > 1$, the tied entries
at rate $1$, and the unaltered entries are unchanged.
Note that $\delta$ is a nonnegative combination of $k$-hot indicator vectors:
\BEQ\label{e-delta}
\delta(n_u,n_t,k)=\frac{1}{\binom{n_t-1}{k-n_u-1}}\sum_{\delta\in\mathcal{I}(n_u,n_t,k)}\delta,
\EEQ
where
\[
\mathcal{I}(n_u,n_t,k) = \left\{\delta\in\{0,1\}^m ~\middle|~ \begin{array}{ll}
    \delta_i=1 & i \leq n_u\\
    \delta_i=0 & i > n_u+n_t\\
    \sum_{i=1}^m \delta_i = k
\end{array}
\right\}
\]
is the set of $k$-hot indicator vectors whose first $n_u$ entries are $1$,
with the remaining $k-n_u$ ones among entries $n_u+1$ through $n_u+n_t$.

\paragraph{Scaling the step.}
We now describe how to compute $s_0$ for the general step.
Let $z(s) = z - s\delta(n_u,n_t,k)$.
The scaling $s_0$ is the minimum of three candidate step sizes,
corresponding to the three events described above.

\begin{enumerate}[label=(\arabic*)]
\item \emph{Merge.}
The step $s_1$ at which the last untied entry reaches the tied value,
\ie, $z(s_1)_{n_u} = z(s_1)_{n_u+1}$.
If $n_u = 0$, there are no untied entries, and $s_1 = \infty$. Otherwise,
$s_1 = \left(a_t - a_u\right)/\left(1-\frac{n_t}{k-n_u}\right)$.

\item \emph{Absorb.}
The step $s_2$ at which the tied value reaches the first unaltered entry,
\ie, $z(s_2)_{n_u+n_t} = z(s_2)_{n_u+n_t+1}$.
If $n_u+n_t = m$, there are no unaltered entries, and $s_2 = \infty$.
Otherwise, since the tied entries decrease at unit rate while $a_e$ is
unchanged, $s_2 = a_t - a_e$.

\item \emph{Termination.}
The step $s_3$ at which the sum of the $k$ largest entries equals $d$,
\ie, $s_3 = (S - d)/w$, where
$w = n_u\frac{n_t}{k-n_u}+(k-n_u)$.
\end{enumerate}
The step is $s_0 = \min\{s_1,s_2,s_3\}$, and the update is
$z \leftarrow z - s_0\delta(n_u,n_t,k)$.

\paragraph{Updating the state.}
Instead of explicitly carrying out this
vector subtraction at every step (which would have a complexity of $O(m)$ per
iterate), the constant size algorithm state is updated to reflect the change
in the vector. The updates are applied in the order shown:
\[
\begin{array}{ll}
    j &\leftarrow j+1\\
    \eta &\leftarrow \eta + s_0\frac{n_t}{k-n_u}\\
    S &\leftarrow S - s_0\left(\frac{n_t}{k-n_u}n_u + (k-n_u)\right)\\
    a_t &\leftarrow a_t - s_0\\
    n_u &\leftarrow \max\{n_u - 1, 0\} \text{~if~} s_0=s_1 \text{~else~}n_u\\
    n_t &\leftarrow \min\{n_t+1,m\}\\
    a_u &\leftarrow z_{n_u} - \eta \text{~if~} n_u > 0 \text{~else~} a_u\\
    a_e &\leftarrow z_{n_u+n_t+1} \text{~if~} n_u+n_t < m \text{~else~} a_e
\end{array}
\]
The first group ($\eta$, $S$, $a_t$) uses the pre-update values
of $n_u$ and $n_t$; the second group ($a_u$, $a_e$) uses the
updated values of $n_u$, $n_t$, and $\eta$.

\subsection{Complexity}
At each step, $n_u$ either decreases by $1$, or remains the same. At each step,
$n_t$ increases by $1$. If the algorithm hasn't terminated after $m$ steps, then
there are $m$ tied entries, thus the subsequent decrease step will reduce the sum
to the desired value, and the algorithm will terminate. Each step of the
algorithm has constant complexity, and thus the algorithm (with sorted input)
has complexity $O(m)$. Therefore, including the sorting and unsorting, the total
complexity of the algorithm is $O(m\log m)$.

\subsection{Correctness}\label{s-correctness}
At termination, the algorithm returns $z = v'-\sum_t\delta_j$.
Note that by \eqref{e-delta}, each $\delta_j$ is a nonnegative combination of
largest-$k$ indicator vectors.
Thus, the algorithm returns $z$ which satisfies
\BEQ\label{e-proj-cond-1}
    v'-z = \sum_{i=1}^{\binom{m}{k}} \lambda_i a_i, \qquad \lambda_i \geq 0, \quad i=1,\ldots,\binom{m}{k}.
\EEQ
Since the decrease step maintains the monotonicity of $v'$,
all previous largest-$k$ indicators are still valid indicators of the largest
$k$ entries of $z$. Thus, the termination criterion provides that
\BEQ\label{e-proj-cond-2}
    a_i^Tz \leq d, \quad \left(a_i^Tz < d \implies \lambda_i = 0\right), \qquad i=1,\ldots,\binom{m}{k}.
\EEQ
However, \eqref{e-proj-cond-1} and \eqref{e-proj-cond-2} are exactly
the KKT conditions for the projection problem described in \S\ref{s-kkt}.
Therefore, the algorithm correctly computes the projection of $v'$ onto the set $\{z\mid f_k(z)\leq d\}$.

\section{Numerical experiments}
\label{s-experiments}
We compare our CVaR projection algorithm and our CVQP solver with two
general-purpose solvers, Mosek and Clarabel. We first benchmark the projection
algorithm on synthetic instances, then benchmark the CVQP solver on
portfolio optimization and quantile regression problems.
Code is available at \url{https://github.com/cvxgrp/cvqp}.

\subsection{Experimental setup}
All experiments were run on a Google Cloud n1-highmem-32 instance
(32 vCPUs, 208 GB RAM), with a two-hour time limit per solve.
For our CVQP solver we use ADMM penalty parameter $\rho^{(0)} = 10^{-2}$,
over-relaxation parameter $\alpha = 1.7$, and update $\rho$ adaptively
with parameters $\mu = 10$ and $\tau = 2$.
We stop when primal and dual residuals satisfy absolute tolerance $10^{-4}$
and relative tolerance $10^{-3}$.
Mosek and Clarabel are called through CVXPY with default settings.
When we tighten our tolerances to $10^{-6}$, objective values agree with both solvers
to about four significant figures.

\subsection{CVaR projection}
We generate random CVaR projection problems with vectors $v\in\reals^m$ having
entries uniformly distributed on $[0,1]$, with $k = (1-\beta)m$ and $\beta = 0.95$.
For each $v$, we set the constraint
threshold to $d = \eta f_k(v)$ in problem \eqref{e-sumklargest_projection},
where $\eta \in (0,1)$ controls problem difficulty (larger $\eta$ gives easier
problems), and project $v$ onto the convex set $\{z\mid f_k(z)\leq d\}$.

\begin{figure}
    \centering
    \includegraphics[width=0.75\textwidth]{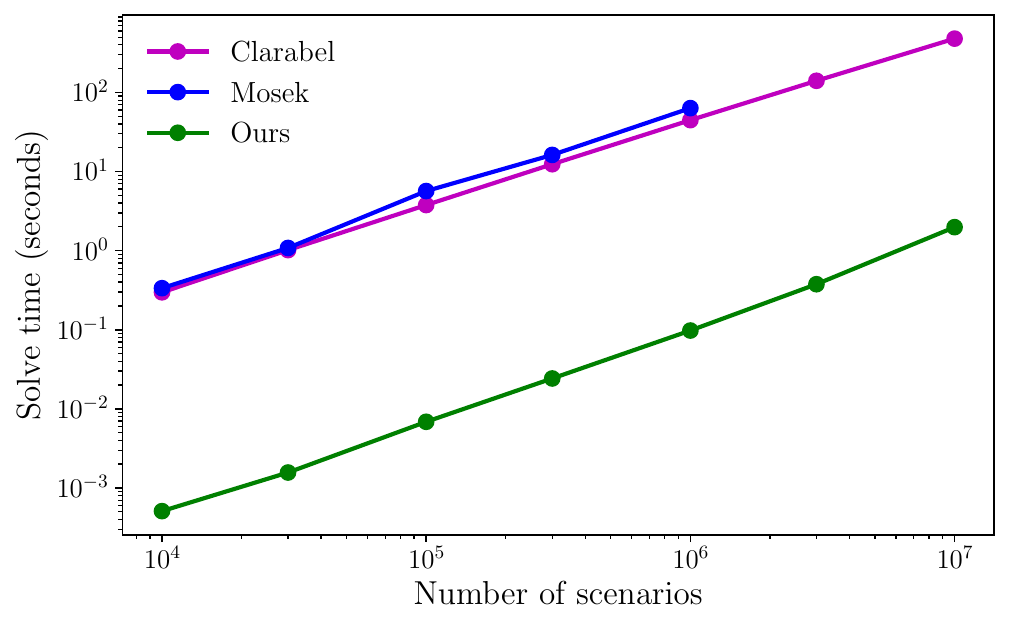}
    \caption{CVaR projection solve times, with $\eta=0.5$.}
    \label{f-cvar-time}
\end{figure}

\paragraph{Results.}
Figure~\ref{f-cvar-time} shows solve times versus number of scenarios $m$ for
$\eta=0.5$, averaged over 10 random instances. Our method solves these problems
in $0.51$ ms at $m=10^4$ and $1.98$ s at $m=10^7$, more than two orders of
magnitude faster than both Mosek and Clarabel across all tested sizes.
Mosek does not solve instances with $m \geq 3\times 10^6$ within the time limit.

\subsection{Portfolio optimization}
We consider a portfolio optimization problem with $n$ assets and $m$ return
scenarios. The matrix $R \in \reals^{m \times n}$ contains asset returns, where
$R_{ij}$ is the return of asset $j$ in scenario $i$. The mean return vector $\mu \in \reals^n$ and covariance matrix
$\Sigma$ are given by
\[
\mu_j = \frac{1}{m}\sum_{i=1}^m R_{ij}, \quad \Sigma = \frac{1}{m}\sum_{i=1}^m (R_i-\mu)(R_i-\mu)^T.
\]

The portfolio optimization problem is
\[
\begin{array}{ll}
\mbox{minimize}   & -\mu^T x + \frac{\gamma}{2}x^T\Sigma x\\
\mbox{subject to} & \mathbf{1}^T x = 1 \\
               & x \geq 0 \\
               & \phi_{\beta}(-Rx) \leq \kappa,
\end{array}
\]
with variable $x \in \reals^n$ (the portfolio weights), risk aversion
parameter $\gamma > 0$, and a CVaR constraint bounding the expected loss
in the worst $(1-\beta)$ fraction of scenarios.

\paragraph{CVQP form.} The problem maps to CVQP form as
\[
\begin{aligned}
P &= \gamma\Sigma, & q &= -\mu, & A &= -R, \\[1em]
B &= \begin{bmatrix} \mathbf{1}^T \\ I \end{bmatrix}, & 
l &= \begin{bmatrix} 1 \\ 0 \end{bmatrix}, & 
u &= \begin{bmatrix} 1 \\ \infty \end{bmatrix}.
\end{aligned}
\]

\paragraph{Problem instances.} We generate return scenarios using a
two-component Gaussian mixture model
\[
R_i \sim \omega \mathcal{N}(\nu \ones, I) + (1-\omega) \mathcal{N}(-\nu \ones, \sigma^2 I), \quad i = 1,\ldots,m,
\]
where $\omega$ is the probability of normal market conditions, $\nu$ is the mean
return per asset (positive in normal conditions, negative in stress periods), and
$\sigma > 1$ scales the volatility during stress periods. We use parameters
\[
\begin{aligned}
\omega &= 0.8, & \nu &= 0.2, & \gamma &= 1, \\
\sigma &= 2, & \beta &= 0.95, & \kappa &= 0.3.
\end{aligned}
\]

\begin{figure}
    \centering
    \includegraphics[width=0.75\textwidth]{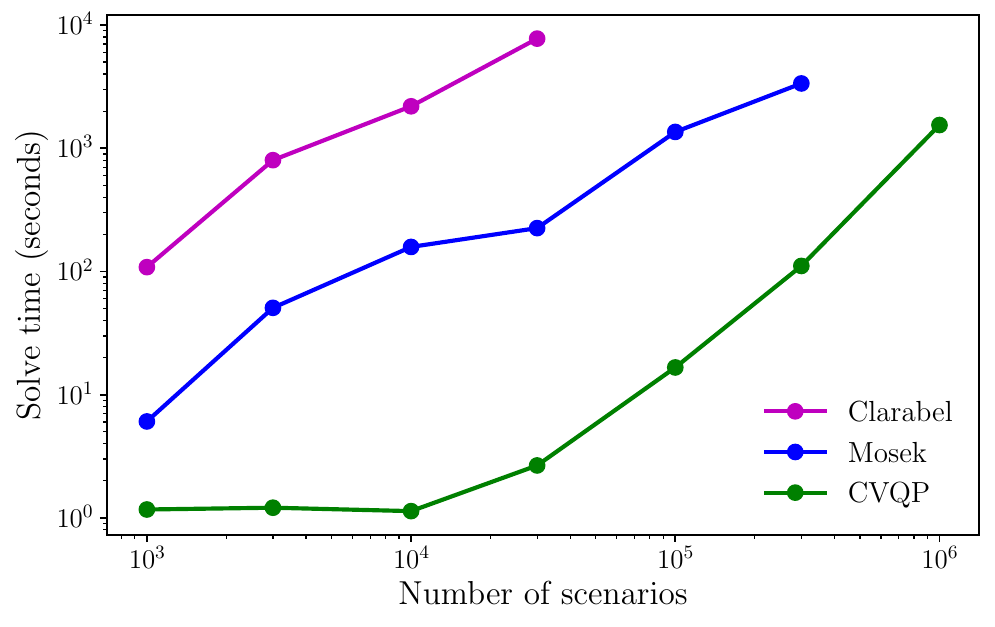}
    \caption{Portfolio optimization solve times, with $n=2000$ assets.}
    \label{f-portfolio-time}
\end{figure}

\paragraph{Results.}
Figure~\ref{f-portfolio-time} shows solve times versus number of scenarios $m$
for problems with $n=2000$ assets, averaged over 3 random instances. Our method
is faster than both general-purpose solvers by one to three orders of magnitude.
Mosek does not solve instances with $m \geq 10^6$ within the time limit,
and Clarabel does not solve instances with $m \geq 10^5$.
Our method solves instances with up to one million scenarios in under 26 minutes.

\subsection{Quantile regression}
Given $m$ data points with features $u_i \in \reals^n$ and responses
$y_i \in \reals$, the $\tau$-quantile regression problem is
\[
\mbox{minimize} \quad \frac{1}{m}\sum_{i=1}^m \rho_\tau(y_i - x^T u_i - x_0),
\]
with variables $x \in \reals^n$ and $x_0 \in \reals$, where
$\rho_\tau(z) = \tau (z)_+ + (1-\tau)(z)_-$ is the tilted $\ell_1$ penalty
and $\tau \in (0,1)$ is the quantile level. The asymmetric penalty weighs
underpredictions $\tau/(1-\tau)$ times more than overpredictions, so the optimal
prediction is the conditional $\tau$-quantile of $y$ given $u$, rather than
the conditional mean.

\paragraph{CVaR reformulation.}
Defining the feature matrix $U \in \reals^{m \times n}$ with rows $u_i^T$,
the response vector $y \in \reals^m$, and $\bar{u} = \frac{1}{m}\sum_{i=1}^m u_i$,
we use the identity $\rho_\tau(z) = (1-\tau)(-z) + (z)_+$ to write the objective as
\[
(1-\tau)\left(x^T\bar{u} + x_0 - \bar{y}\right) + \frac{1}{m}\sum_{i=1}^m (y_i - u_i^Tx - x_0)_+,
\]
where $\bar{y} = \frac{1}{m}\sum_{i=1}^m y_i$. The terms involving $x_0$ can be written as
\[
(1-\tau)\left[x_0 + \frac{1}{(1-\tau)m}\sum_{i=1}^m(y_i - u_i^Tx - x_0)_+\right],
\]
which is $(1-\tau)$ times the sample CVaR formula evaluated at $\alpha = x_0$.
Minimizing over $x_0$ yields the infimum in the CVaR definition.
Since $\bar{y}$ is constant and $(1-\tau) > 0$, the problem reduces to
\[
\mbox{minimize} \quad x^T \bar{u} + \phi_\tau(-Ux + y),
\]
with variable $x \in \reals^n$. Introducing an epigraph variable $t$ for the
CVaR term gives
\[
\begin{array}{ll}
\mbox{minimize}   & x^T \bar{u} + t \\
\mbox{subject to} & \phi_\tau(-Ux + y) \leq t.
\end{array}
\]

\paragraph{CVQP form.}
We introduce an auxiliary variable $s$ with the constraint $s=1$, and define
$\tilde{x} = (x, t, s) \in \reals^{n+2}$. The CVQP parameters, with CVaR
level $\beta = \tau$, are
\[
\begin{aligned}
P &= 0, & q &= \begin{bmatrix} \bar{u} \\ 1 \\ 0 \end{bmatrix}, & A &= \begin{bmatrix} -U & -\ones & y \end{bmatrix}, \\[1em]
B &= \begin{bmatrix} 0^T & 0 & 1 \end{bmatrix}, &
l &= 1, &
u &= 1.
\end{aligned}
\]

\paragraph{Problem instances.}
We generate features $u_i \in \reals^n$ with independent standard normal
entries, a vector with entries $\beta_j \sim \mathcal{N}(0, 1/(1+j))$,
and responses $y_i = u_i^T \beta + \epsilon_i$, where $\epsilon_i$ follows a
$t$-distribution with 5 degrees of freedom scaled by $0.1$. We use quantile
level $\tau = 0.9$.

\begin{figure}
    \centering
    \includegraphics[width=0.75\textwidth]{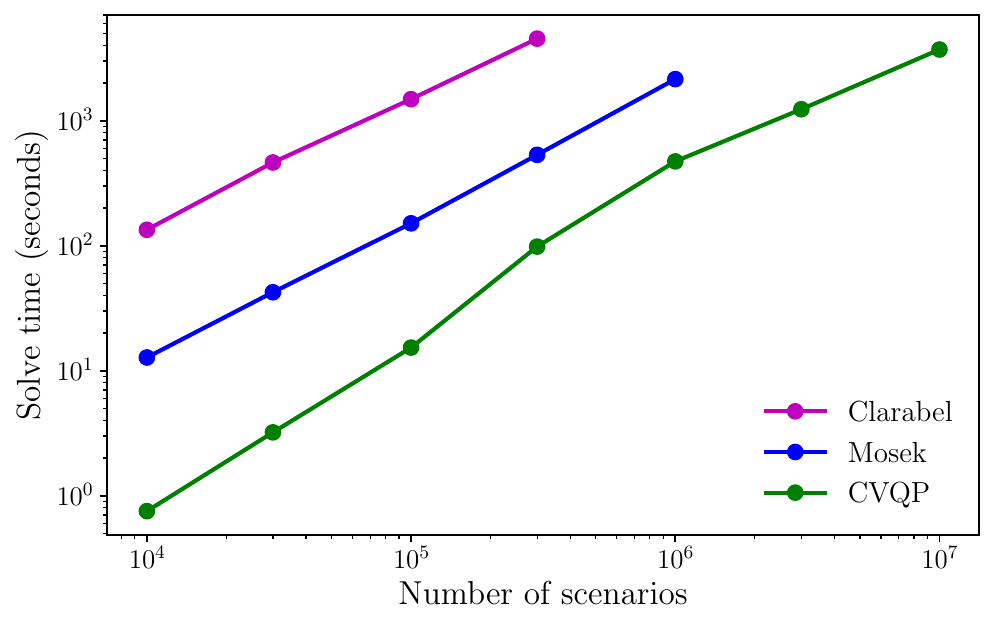}
    \caption{Quantile regression solve times, with $n=500$ features.}
    \label{f-qr-time}
\end{figure}

\paragraph{Results.}
Figure~\ref{f-qr-time} shows solve times versus number of scenarios $m$ for
quantile regression problems with $n=500$ features, averaged over 3 random
instances. Our method is faster than both Mosek and Clarabel across all tested sizes,
by a factor ranging from about 5 to over 100.
Mosek does not solve instances with $m > 10^6$ within the time limit,
and Clarabel does not solve instances with $m \geq 3\times 10^5$.
Our method scales to $m=10^7$ scenarios.

\section{Extensions and variations}
\label{s-extensions}
Our method extends to more general problems and admits several practical variations.

\paragraph{More general constraints on $\tilde{z}$.}
When projection onto a set $\mathcal{C}$ is efficient (via a closed-form solution
or a fast algorithm), we can directly handle constraints of the form
$\tilde{z} \in \mathcal{C}$ in place of the box constraints $l \leq \tilde{z} \leq u$.

\paragraph{Extensions to the $x$-update.}
The method extends to additional constraints on $x$ or
nonquadratic objective functions. The $x$-update then becomes
\[
\begin{array}{ll}
\mbox{minimize} & f(x) + (\rho/2)\|Ax - z^k + u^k\|^2_2 + (\rho/2)\|Bx - \tilde{z}^k + \tilde{u}^k\|^2_2\\
\mbox{subject to} & Gx \in \mathcal{D},
\end{array}
\]
where $f$ is any convex function (not just quadratic), $G \in \reals^{p'
\times n}$, and $\mathcal{D} \subseteq \reals^{p'}$ is a convex set. Using a
standard solver with warm-starting, the per-iteration cost remains manageable, 
though typically higher than solving the original linear system.  

\paragraph{Equilibration.}
Equilibration as a preprocessing step often improves
ADMM convergence, particularly for problems with poorly scaled data.

\paragraph{Warm starting.}
Two effective strategies for choosing the initial point are
(1) solving with a reduced set of scenarios and using that solution to initialize the full problem,
and (2) using a quadratic approximation of the CVaR constraint to compute an initial point.
Similarly, the sorting step in the CVaR projection can be warm-started
using the sorted order from the previous ADMM iteration.
Since the projection input changes only slightly between iterations, the
sorted order is nearly preserved, and an adaptive sorting algorithm can
exploit this to reduce the per-iteration cost in practice.

\section*{Acknowledgments}
David Pérez-Piñeiro received support from the Research Council of Norway and the
HighEFF Research Centre (Centre for Environmentally Friendly Energy Research,
Grant No. 257632).

\clearpage
\printbibliography

\end{document}